\documentclass[11pt]{amsart}

\usepackage{amsmath,amssymb,amsfonts,amsthm}

\newcommand{\N}{\ensuremath{\mathbb{N}}}

\newcommand{\set}[1]{\ensuremath{\left\{#1\right\}}}

\newtheorem{thm}{Theorem}

\newtheorem{cor}{Corollary}

\usepackage[colorlinks=False,linkcolor=blue]{hyperref}
\usepackage[utf8]{inputenc}
\usepackage{enumerate}
\usepackage{xcolor}
\usepackage{url,cite}
\usepackage{mathtools,minibox}

\newcommand{\upc}[1]{{\MakeUppercase #1}}

\newcommand{\st}{\ensuremath{\colon}}

\title{New developments of the Odds Theorem}
\author{Rémi Dendievel}
\address{Université Libre de Bruxelles\\
Campus Plaine, CP 210\\
B-1050 Bruxelles\\
Belgique}

\begin{document}


\keywords{Odds algorithm, records, secretary problems, robotic
  maintenance, clinical trials, investment problems, multiplicative
  odds, Markov chains.}
\subjclass{60G40}

\maketitle

\vspace{-5ex}
\begin{center}
  Université Libre de Bruxelles\\
  Faculté des Sciences\\
  Département de Mathématique
\end{center}

\begin{abstract}
  The odds theorem and the corresponding solution algorithm (odds
  algorithm) are tools to solve a wide range of optimal stopping
  problems. Its generality and tractability have caught much
  attention. (Google for instance ``Bruss odds'' to obtain a quick
  overview.) Many extensions and modifications of the this result have
  appeared since publication in~2000. This article reviews the
  important gnew developments and applications in this field.  The
  spectrum of application comprises as different fields as secretary
  problems, more general stopping problems, robotic maintenance
  problems, compassionate use clinical trials and others.

  This review also includes a new contribution of our own.
\end{abstract}

\section{The Original Odds Theorem}
\label{sec:orig-odds-theor}

The odds-theorem is a result in the theory of optimal stopping which
can be applied for many interesting sequential decision problems. The
original version of the Odds-algorithm is the work of Bruss~(2000). He
discovered it when he saw common features between quite different
looking optimal stopping problems.

The framework is as follows. There are $n$ random variables which are
observed sequentially. It is desired to stop online with maximum
probability on a last specific event. No recall is permitted. Here
``specific'' is understood as being defined in terms of the interest
of the decision maker. Such problems can be readily translated into
the equivalent problem of stopping on a sequence of independent
indicators. The parameters of those indicators variables are supposed
to be known. Maximizing the resulting objective function means then
maximizing the probability of selecting the last indicator of a
specific event to be equal to~1 from the sequence. For convenience a
variable equal to one will be called a \emph{success}.

The independence property can be relaxed but is, at least locally,
important. This is the main reason why the optimal strategy turns out
to be a threshold rule based on a fixed time index.

We first recall the odds theorem, upon which the Odds-algorithm is
based.
\begin{thm}[Odds-Theorem, Bruss (2000)]\label{thm:odds-theorem}
  Let $I_1,I_2,\dotsc,I_n$ be $n$ independent Bernoulli random
  variables, with $n$ known. We denote ($i=1,2,\dotsc,n$) $p_i$ the
  parameter of the variable $I_i$ ($p_i\in[0,1]$). Let $q_i=1-p_i$ and
  $r_i=p_i/q_i$. Define the index
  \begin{displaymath}
    s = \begin{cases} \max\set{k\in\set{1,2,\dotsc,n}\st \sum_{j=k}^n r_k
        \geq 1} & \text{ if $\sum_{i=1}^n r_i\geq 1$,} \\
      1 & \text{ otherwise.}
    \end{cases}
  \end{displaymath}
  To maximize the probability of stopping on the last ``1'' of the
  sequence it is optimal to stop on the first ``1'' that we meet among
  the variables $I_s,I_{s+1},\dotsc,I_n$.
\end{thm}

Based on the this theorem, we can describe the \emph{odds-algorithm} as follows:
\begin{quote}
  \begin{enumerate}[1{$^\circ$}]
  \item Compute the odds $r_j$ for $j=n,n-1,\dotsc$ successively;
    compute the threshorld $s$ easily by looking at the running sum
    $r_n+r_{n-1}+\dotsb+r_j$ and stop as soon as this sum reaches or
    exceeds $1$. This defines the stopping threshold $s$ with $n-s+1$
    further variables to come.  We then must wait for the first $k\ge
    s$ (if any) with $I_k=1$. Otherwise we must stop at index $n$
    anyway.
  \item The optimal win probability is given by $V(n)=R_sQ_s$, where
    $Q_s=\prod_{j=s}^n q_j$ and $R_s = \sum_{j=s}^n r_j$.
 \end{enumerate}%
\end{quote}%
A subsequent article of Bruss~(2003) gives lower and upper bounds for
the quantity $V(n)$.

The algorithm provides thus the optimal strategy $\tau_s$ and is
optimal with respect to other considerations: linear complexity, ease
of computation. Indeed with simple values of $p_k$ this can even be
computed by head.

\subsection*{Examples of applications}
\label{sec:exampl-appl}

This result has been immediately noticed for its simplicity and
generality. Several generalizations have appeared since then. The
overview of these works is the content of the next sections. We recall
first a few applications included in Bruss~(2000).

\subsubsection*{1) Classical Secretary Problem}

An interviewer has an appointment with each of the $n$ (fixed)
secretaries who are applying for a certain job. The secretaries'
quality are independent of each other. He is not able to quantify the
quality of a secretary but he is able to rank them relatively to each
other. That is, after having observed, say, $k$ secretaries one after
the other he can compute their relative ranks, and in particular he is
able to remember which one was the best candidate among these first
$k$ secretaries.

Solving this problem with the odds theorem is straightforward. An
observation at time $k$ is a record if and only if his relative rank
among the first $k$ variables is 1. It is known that $p_k=1/k$. Hence
$q_k=(k-1)/k$ and $r_k=1/(1-k)$. The $s$ index in the odds theorem can
be computed and gives the optimal strategy. It is the largest $s$ such
that $\sum_{j=s}^n (j-1)^{-1} \ge 1$. Therefore $V(n)=\frac{s-1}n
\sum_{j=s-1}^{n-1}j^{-1}$. Note that $s=s(n)\sim n / e$ and thus
$V(n)\to 1/e$. Indeed $V(n)\ge 1/e$ for all $n=1,2,\dots$.

We should also mention here that even if the parameter $n$ is unknown,
the best candidate can always be selected with probability of at least
$1/e$, and this is a very tractable model. This is the so-called
$1/e$-law of Bruss~(1984).

For a related problem see also Suchwa\l ko and Szajowski~(2003).

\medskip
\noindent%
\textbf{Remark.} In the setting of the odds theorem, the real number
$1/e$ is a lower bound for $V(n)$ provided that $\sum_{i=1}^n r_i \ge
1$ (see Bruss (2003)). This condition is always met in the classical
secretary problem.

\subsubsection*{2) Grouping observations}

Hsiao and Yang (2000) adress and solve secretary problem with group
interviews. The selection is considered a success if the selected
group that contains the best of all observations. This can also
immediately be solved using the odds algorithm.

\subsubsection*{3) Refusal probability}

It is easy to introduce an availability probability in the
framework. If the probability for some variable $I_k$ to be equal to 1
is $p_k$ and, independently of the value of $I_k$, the variable is
available with probability $a_k$, then the probability of having an
available 1 is $\tilde p_k = p_k a_k$. The odds algorithm can compute
the strategy based on the $\tilde p_k$'s.

\subsubsection*{4) Unknown number of observations}

This model can immediately be generalized to a random number
$N$ of events as long as we assume independance of the indicators of
the successes. For instance we model the unknown number $N$ by a
time-embedding through
\[
p_k=P(I_k=1|\exists\text{\small{} an observation at time $k$})\cdot
P(\exists \text{\small{} an observation at time $k$})\,.
\]


\subsubsection*{Dice game}

A well-known game is the following. A die is thrown $N$ times, $N$
fixed. To one player it is asked to bet an amount of money one of the
$N$ throws. He wins if the die shows the value 6 at that time and if
there is no more 6's in the following throws.

Since the probability of obtaining a 6 at any time is $1/6$, we input
the values $(1/6,1/6,1/6,\dots)$ as the parameteres $(p_k)$ in the
algorithm. We obtain $s=N-4$. That is, we will only look at the value
shown by the last 5 dice and bet our money as soon as we see a 6.

\medskip

Additional applications will be outlined in the last section.

\section{Stopping On The Last Success: Unkown Odds, Random Length,
  Recall and Other Modifications}
\label{sec:stopp-last-succ}

\subsection{Unkown odds}
\label{sec:unkown-odds}

This section addresses a difficult problem. What would we know if all
we knew was that the observed variables are independent but have an
unknown parameter? A natural approach would be to estimate the odds
sequentially and to plug the estimates into the odds algorithm. Let us
call such a strategy an ``empirical odds strategy''.

For a detailed study of the performance of empirical odds strategies
we refer to the paper of Bruss and Louchard (2009) in which they
analyze and test several modifications.

The intuition tells us that the optimal strategy lies in the class of
empirical strategies. But at the moment there is no theoretical result
to show this. In a simple setup (small $n$) dynamic programming shows
that this intuition is correct that is, that the empirical odds
strategy is optimal. However for a larger $n$ this is still an open
problem, and an important one regarding applications, for it is often
closer to reality than the model in which we assume that we know the
parameters.

For motivation and examples of application for this model, see
also Bruss (2006). In particular this also treats important
applications in the domain of clinical trials.

\subsection{Stopping on the $m$-th last success}
\label{sec:stopping-m-th}

The paper of Bruss and Paindaveine~(2000) follows the spirit of the
original odds-theorem of~2000. The setting is the same as in the
previous section. The objective is now to predict the $m$-th last
success upon its arrival, that is, to find the stopping time $\tau$
that maximizes the probability
\begin{equation}
  \label{eq:7} P\left( \sum_{k=\tau}^{n} I_k = m \right)
\end{equation} and the optimal strategy associated with this win
probability.


\begin{thm}[Bruss and Paindaveine (2000)] An optimal rule for stopping
  on the $m$-th last success exists and is to stop on the first index
  (if any) with $I_k=1$ and $k\ge s_m$ for some fixed
  $s_m\in\set{1,2,\dotsc,n-m+1}$.
\end{thm}
The $s_m$ are computed as follows: define
\begin{align*} \pi_k & \coloneqq \# \set{j\ge k | r_j > 0}\\
\intertext{and}R_j^{(k)} &= \sum_{k\leq i_1<i_2<\dotsb<i_j\leq n}
r_{i_1}r_{i_2}\dotsb r_{i_j}\,.
\end{align*}

We then have
\begin{equation}
  \label{eq:8} s_m = \sup\set{1, \sup\set{1\le k\le n-m+1 : R_m^{(k)}
\ge mR_{m-1}^k \text{ and } \pi_k \ge m} }.
\end{equation}

It is to mention that unimodality of the optimal strategy is
not straightforward and that it needs a delicate treatment.  The
stopping index $s_m$ can be computed but has a slightly more
sophisticated form than before.

\subsection{Hsiao and Yang's Markovian framework}
\label{sec:hsiao-yangs-mark}

\subsubsection*{Homogeneous case}
\label{sec:hoomgeneous-case}

Hsiao and Yang~(2002) study a modification of the same model where now
the $I_1, I_{2},\dotsc, I_N$ form a Markov chain. The authors prefer
to renumerate the indicators backwards. Hence let $I_N,I_{N-1}, \dots,
I_{1},I_0$ be a Markov chain with the following structure:
\begin{gather*} P(I_{n-1}=1|I_{n}=0) = \alpha_{n}\\
P(I_{n-1}=0|I_{n}=1) = \beta_{n}
\end{gather*} The authors treat again the objective to stop with
maximum probability on the last success. Let
\begin{gather*} S_n = P(I_{j}=0, \text{ $ \forall
j=n-1,\dots,0$}|I_n=1) = \beta_n\prod_{i=1}^{n-1}(1-\alpha_{i})\\
\intertext{and let} q_0(n) = \text{optimal success probability on
$I_{n-1},\dots,I_0$ given that $I_n=0$},\\ q_1(n) = \text{optimal
success probability on $I_{n-1},\dots,I_0$ given that $I_n=1$}.
\end{gather*}

It can be seen that the stopping time is defined as the first $n$ such
that
\begin{equation} 
  S_n \ge q_1(n)\,,\label{eq:1}
\end{equation}
and as $0$ if there is no such $n$. Let $(\phi_j,j=N,\dotsc,1,0)$ be
the stopping strategy. So this sequence is adapted to the process
$(I_j, j=N,\dots,1,0)$; $\phi_j=1$ means that we choose to stop at
time $j$ if $I_j=1$ and $\phi_j=0$ means that we continue observing
more variables, whatever the value of $I_j$. We always have $\phi_0=1$
because in our problems a decision must be made within the set
$\{N,N-1,\dotsc,0\}$.

The first result is obtained in the case where $\alpha_n$ and
$\beta_n$ are constants, for all $n$; we set $\alpha:=\alpha_0$ and
$\beta:=\beta_0$.
\begin{thm}[Hsiao and Yang (2002)] If $\beta
  \in [\frac12, 1]$, then
  \begin{enumerate}[{\upshape (i)}]
  \item if $\alpha=0$, $\phi_j=1$ for all $j=N,N-1,\dots,0$;
  \item if $\alpha=1$, $\phi_0=\phi_1=1$ and $\phi_j=0$ for
    $j=N,N_1,\dots,2$;
  \item if $\alpha\in(0,1)$, $\phi_j=0$ for
    $j\in\set{N,\dots,r+2,r+1}$ and $\phi_j=1$ for
    $j\in\set{r,r-1,\dots,0}$, where $r= \min\set{\lfloor
      (\beta-2\alpha)(1-\alpha)/\alpha\beta\rfloor+2, N}$
  \end{enumerate} Therefore, there exists an $r$ such that $\tau_N =
  \sup\set{0 \le i \le N | I_i = 1, i \le r}$ with the convention that
  $\sup \varnothing = 0$.
\end{thm}

The case which involves more calculations is the third. Hsiao and Yang
obtain the explicit form of $q_1(k)$ and $q_0(k)$ by solving the
recurrence defining those two functions. The index $r$ is then
obtained using~(\ref{eq:1}) and replacing $q_1$ by its explicit value.

We should explain why this result holds for $\beta\in[\frac12, 1]$. A
high value for $\beta$ means that once we observe the variable with
value 1, it is likely that there will be other variables in the future
equal to 0. On the other hand, a small $\beta$ means that the it is
likely that after a 1, we have many variables being equal to 1.

Similar but more delicate cases arise for the case
$\beta\in(0,\frac12)$. The resulting strategy in the non-degenerate
case deserves also an explanation. Depending on whether the quantity
\begin{equation} (\alpha+\beta)\beta(1-\alpha)^n -
\beta(1-\alpha-\beta)^{n+1}\label{eq:2}
\end{equation} is smaller than $\alpha$ for any $n=0,1,\dots,N-1$ or
exceeds $\alpha$ for some $n$, the resulting strategy has a very
different structure. We now state their result.

\begin{thm}[Hsiao and Yang (2002)]Let $\beta\in(0,\frac12)$.\hfill\phantom{,}
  \begin{enumerate}[{\upshape (i)}]
  \item If the quantity defined in {\upshape(\ref{eq:2})} is always
    smaller than $\alpha$, we have $\phi_0=1$, $\phi_j=0$ for all $j>0$;
  \item If $\alpha\neq 0$ and there exists an $r\in\set{N-1,\dots,1}$
    such that
    \begin{equation} 
      (\alpha+\beta)\beta(1-\alpha)^r -
      \beta(1-\alpha-\beta)^{r+1} \ge \alpha >
      (\alpha+\beta)\beta(1-\alpha)^k -
      \beta(1-\alpha-\beta)^{k+1}\label{eq:3}
    \end{equation} for all $k<r$. Define
    \begin{displaymath} m = \left\lfloor
        \frac{(\alpha+\beta)(\alpha^2-\alpha+\beta)(1-\alpha)^{r-1}-
          a[1-(1-\alpha-\beta)^{r+1}]}{\alpha\beta(\alpha+\beta)
          (1-\alpha)^{r-1}}\right\rfloor + 1,
    \end{displaymath} we have the following optimal strategy
    \begin{equation*}
      \begin{gathered} \phi_j=
        \begin{cases} 1 & \text{for $j \in \set{0,
              r+1,r+2,\dots,r+m}$}\\ 0 & \text{else}
        \end{cases} \hfill\text{ if $r+m<N$,} \\ \phi_j=
        \begin{cases} 1 & \text{for $j \in \set{0,
              r+1,r+2,\dots,N}$}\\ 0 & \text{else}
        \end{cases} \hfill\text{ if $r+m\ge N$.}
      \end{gathered}
    \end{equation*}
    This is a case where the optimal stopping has more than one
    stopping island. 
  \item If {\upshape(\ref{eq:3})} is verified for some $r$ for all
    $k<r$ and if $\alpha=0$, then there exists $r<N$ such that
    \begin{displaymath} \phi_j =
      \begin{cases} 1 & \text{for } j \in\set{0,r,r+1,r+2,\dots,N} \\
        0 & \text{for } j \in\set{1,2,\dots,r-1}
      \end{cases}\,.
    \end{displaymath} This strategy represents the stopping time
    $\tau_N$ defined as
    \begin{displaymath} \tau_N = \sup\set{0\le i \le N | I_i = 1,\,
        i=0 \text{ or } i > r}\,.
    \end{displaymath}
  \end{enumerate}
\end{thm}

The probability of selecting the last success by using the optimal
strategy can be computed for all $\beta$.

\subsubsection*{Nonhomogeneous case}
\label{sec:nonhomogeneous-case}

Hsiao and Yang then study the corresponding non-homogeneous case and
obtain the following theorem under some assumptions.

\begin{thm}[Hsiao and Yang (2002)] If $\alpha_n+\beta_n\ge 1$ for all
  $n$ then
  \begin{gather*} \phi_j =
    \begin{cases} 1 & \text{for } j \in\set{0,1,\dots,r}\\ 0 &
\text{for } j \in\set{r+1, r+2,\dots,N}
    \end{cases}
  \end{gather*} where
  \begin{displaymath} r = \inf\biggl\{k\in\mathbb N_N\cup\{0\}:
    \sum_{l=1}^k
    \frac{\alpha_l\beta_{l-1}}{(1-\alpha_l)(1-\alpha_{l-1})}
    +\frac{\beta_k(1-\beta_{k+1})}{\beta_{k+1}(1-\alpha_k)} > 1
    \biggr\}\,.
  \end{displaymath}

\end{thm}

We can remark that in their results the optimal strategy cannot be
simplified into a ``sum-the-odds'' strategy (any kind of odds).



The optimal strategy may have now a different form. One can see that,
as pointed out in the theorem, there can be more than just one
stopping island.

\subsection{Tamaki's Markovian result}
\label{sec:tamak-mark-result}

Tamaki~(2006) tackles a similar problem as in the previously exposed
Markovian framework of Hsiao and Yang. There are important
differences, however. First the hypotheses on the transition
probabilities are different. Second, his objective is to obtain a
solution relying on a sum of odds.

Let $I_1,I_2,\dotsc,I_n$ be a sequence of independent indicator
variables. Let us study the following Markov dependence between the
variables:
\begin{align*} \alpha_j &= P(I_{j+1}=1|I_j=0)\\ \beta_j &=
P(I_{j+1}=0|I_j=1)
\end{align*} for $1 \leq j\leq n-1$.  And let us assume that
$\alpha_n=0$, $\beta_n=1$.  We write $\bar\alpha_j$ and $\bar\beta_j$
for $1-\alpha_j$ and $1-\beta_j$ respectively. The result is as
follows:
\begin{thm}[Tamaki (2006)] Assume that
  \begin{enumerate}[a)]
  \item $(\alpha_j)$ is non-increasing in $j$,
  \item $(\beta_j)$ is non-decreasing concave in $j$.
  \end{enumerate} Then an optimal rule stops on the first index $k\ge
  s$ such that $I_k=1$ and where
  \begin{displaymath} s = \sup \left\{1\leq k\leq n :
\frac{\bar\beta_k}{\beta_k} \frac{\beta_{k+1}}{\bar\alpha_{k+1}} +
\sum_{j=k+1}^{n-1} \frac{\alpha_j}{\bar\alpha_j}
\frac{\beta_{j+1}}{\bar\alpha_{j+1}} \ge 1 \right\}
  \end{displaymath} with the natural convention that the empty sum
equals 0.
\end{thm}

\subsection{Multiple sum-the-odds theorem (Ano et al., 2010)}
\label{sec:multiple-sum-odds}

Suppose that we are given $m\in\N$ selection chances in the problem
described in the preceding section. Let $V_i^{(m)}$ , $i\in\N$ ,
denote the conditional maximum probability of win provided that we
observe $X_i = 1$ and select this success when we have at most $m$
selection chances left. Let $W_i^{(m)},i\in\N$, denote the conditional
maximum probability of win provided that we observe $X_i = 1$ and
ignore this success when we have at most $m$ selection chances
left. Let, furthermore, $M_i^{(m)}$, $i\in\N$, denote the conditional
maximum probability of win provided that we observe $X_i = 1$ and are
faced with a decision to select or not when we have at most $m$
selection chances left. The optimality equation is then given by
\begin{displaymath} M_i^{(m)} = \max\set{V_i^{(m)}, W_i^{(m)}},\quad
  i\in\N.
\end{displaymath}


For each $i\in\N$, define recursively the quantities $H_i^{(m)}$ by
\begin{align}
  H_i^{(1)} &= 1 - \sum_{j=i+1}^N r_j, \\
  H_i^{(m)} &= H_i^{(1)} + \sum_{j=(i+1)\vee i_*^{(m-1)}}^{N} r_j
  H_j^{(m-1)},
\end{align}
where $i_*^{(m)} = \min\set{i\in\N: H_i^{(m)}>0}$.

Now the theorem in Ano et al. reads:
\begin{thm}[Ano, Kakinuma, Miyoshi (2010)] Suppose that we have at
  most $m\in\N$ selection chances. Then, the optimal selection rule
  $\tau_*^{(m)}$ us given by
  \begin{equation}
    \tau_*^{(m)} = \min\set{i\ge i_*^{(m)}:X_i=1}\\
  \end{equation}
  where $\min \varnothing = +\infty$. Furthermore, we have
  \begin{equation}
    \label{eq:13}
    1\le i_*^{(m)} \le i_*^{(m-1)} \le \dotsb \le i*_{(1)} \le N.
  \end{equation}

\end{thm} It would be interesting to have an intuitive understanding
of the quantities $H_i^{(m)}$, but this seems difficult.

In Ano and Matsui (2012), a lower bounds for the multiple stopping
problem is obtained.

\subsection{Random Length}
\label{sec:random-length-tamaki}

Tamaki, Wang and Kurushima (2008) allow random length and provide a
sufficient condition for the optimal rule to be of threshold type.

Ano, in a preprint (2011), tackles again the multiple stopping
problem, with random length.

Random length and refusal probability at the same time are studied in
Horiguchi and Yasuda (2009).

\subsection{Ferguson's modification of the Odds-Theorem}
\label{sec:fergusons-model}

Ferguson~(2008) proposed the following modification of the original
odds-theorem described in section~\ref{sec:orig-odds-theor}.

Let $Z_1,Z_2,\dots$ be a stochastic process on an arbitrary space with
an absorbing state called 0. For $i=1,2,\dotsc$, let $Z_i$ denote the
set of random variables observed after succes $i-1$ up to and
including success $i$. If there are less than $i$ successes, we let
$Z_i=0$, where 0 is a special absorbing state. The general model is as
follows.
 
We make the assumption that with probability one the process will
eventually be absorbed at 0. We observe the process sequentially and
wish to predict one stage in advance when the state 0 will first be
hit. If we predict correctly, we win 1, if we predict incorrectly we
win nothing, and if the process hits~0 before we predict, we win
$\omega$, where $\omega<1$. This is a stopping rule problem in which
stopping at stage $n$ yields the payoff
\begin{equation}\label{eq:4}
  \begin{split} Y_n &= \omega I(Z_n=0) + I(Z_n\neq 0)P(Z_{n+1}=0 |
\mathcal{G}_n)\qquad \text{for $n=1,2,\dots$} \\ Y_\infty &= \omega
  \end{split}
\end{equation} 
where $\mathcal G_n=\sigma(Z_1,\dots,Z_n)$, the $\sigma$-field
generated by $Z_1,\dots,Z_n$. The assignment $Y_\infty=\omega$ means
that if we never stop, we win $\omega$.

The resulting \emph{one-stage look-ahead rule} (1-sla) is to stop at
index $N$ defined by
\begin{equation} N \coloneqq \min\set{k : Z_k =0 \text{ or } (Z_k\neq
0 \text{ and } W_k/V_k \leq 1-\omega)}
\end{equation}
where
\begin{align*} \qquad V_k &= P\set{Z_{k+1}=0 | \mathcal{G}_k}, \\
  W_k &= P\set{Z_{k+1}\neq0, Z_{k+2} = 0 | \mathcal{G}_k }.
\end{align*}
The event $\set{Z_{k+1}\neq0, Z_{k+2} = 0}$ given the history
$\mathcal{G}_k$ describes the event that there is exactly one success
in the future because 0 is an absorbing state. A sufficient condition
for the problem to be monotone is, as Ferguson shows,
\begin{equation}\label{eq:5} W_k/V_k \text{ is a.s. non-increasing in
$k$}.
\end{equation}
\begin{thm}[Ferguson (2008)]\label{sec:ferg-modif-odds} Suppose that
  the process $Z_1,Z_2,\dotsc$
  has an absorbing state 0 such that $P(Z_k \text{ is absorbed at 0})
  = 1$ and that the stopping problem with reward sequence~{\upshape
    (\ref{eq:4})} satisfies the condition~{\upshape
    (\ref{eq:5})}. Then the 1-sla is optimal.
\end{thm}

This model enables us to tackle more general problems. The level of
generality and abstraction of this model makes it a very tractable
result. Furthermore, Ferguson's paper contains several examples for
which the 1-sla rule turns out to be a sum-the-odds strategy.

\section{The Role of the $k$-fold Multiplicative Odds}
\label{sec:the-role-of}

As before in the paper of Bruss and Paindaveine~(2000) where the
authors considered a group of last successes, we have here again to
deal with multiplicative odds. We now present the problem studied by
Tamaki~(2000).


The problem is the following. Find the strategy that maximizes the
probability of stopping on any of the last $m$ successes. All
hypotheses are the same as the ones mentioned in the original odds
theorem in Section~\ref{sec:orig-odds-theor}.

In Ferguson's framework from Section~\ref{sec:fergusons-model}, the
current problem leads us to consider the following payoffs
\begin{displaymath}
  Y_k = I(Z_k\neq 0) P(Z_{k+m} = 0 | \mathcal F_k),\quad
  k=1,2,\dotsc,n,
\end{displaymath}
where $\mathcal F_k$ is the $\sigma$-field generated by
$Z_1,Z_2,\dots,Z_k$, and the following quantities
\begin{align*}
  V_k &= P(Z_{k+1} = 0 | \mathcal F_k), \\
  W_k &= P(Z_{k+m}\neq 0, Z_{k+m+1}=0 |\mathcal F_k).
\end{align*}

A corollary of Theorem~\ref{sec:ferg-modif-odds} from
Section~\ref{sec:fergusons-model} is as follows:
\begin{cor}[Ferguson (2008)]
  Suppose that $n$ Bernoulli random variables $X_1,X_2,\dots,X_n$ are
  observed sequentially. Let $\mathcal F_1,\mathcal F_2,\dots,\mathcal
  F_n$ be an increasing sequence of sigma-fields such that
  $\set{X_j=1}\in\mathcal F_j$ for all $1\le j\le n$. Let
  \begin{align*}
    V_k &= P(X_{k+1} + \dotsb + X_n = 0 | \mathcal F_k),\\
    W_k &= P(X_{k+1} + \dotsb + X_n = m | \mathcal F_k).
  \end{align*}
  Then the optimal rule is determined by the stopping time
  \begin{displaymath}
    N_m = \min\set{k\ge 1 : X_k = 1 \text{ and } W_k/V_k \le 1},
  \end{displaymath}
  provided that the sequence $(W_k/V_k,k\ge 1)$ is monotone
  non-increasing.
\end{cor}

Set, as above, 
\begin{equation}\label{eq:6}
R_j^{(k)} = \sum_{k\leq i_1<i_2<\dotsb<i_j\leq n}
  r_{i_1}r_{i_2}\dotsb r_{i_j}.
\end{equation}
The result is a sum-the-odds strategy, but the odds are the
multiplicative odds given in (\ref{eq:6}).

\begin{thm}[Tamaki (2010)] For the stopping
  problem of maximizing the probability of stopping on any of the last
  $m$ successes in $n$ independent Bernoulli trials, the optimal rule
  stops on the first success $X_k=1$ with $k\ge s_m$, if any, where
  \begin{displaymath}
    s_m = \min\set{k\ge 1 : R_{m}^{(k+1)} \le 1}.
  \end{displaymath}
  Moreover, the maximal probability of win is
  \begin{displaymath}
    v_m = \left( 
      \prod_{s_m}^{n} q_j
    \right)
    \left( 
      \sum_{s_m}^{n} R_{j}^{s_m}
    \right).
  \end{displaymath}
\end{thm}

\noindent\textbf{Remark.} The two problems we mentioned here involve the
$m$ last 1's of the sequence. It would be interesting to know in
advance if a particular problem will be a sum-the-odds theorem
involving the multiplicative odds.

\medskip

We now include a new contribution.

\section{A non-informative problem in continuous time}
\label{sec:non-inform-probl}

It is possible to translate the problem described in
section~\ref{sec:unkown-odds} to a continuous setting as follows. Let
$I_1,I_2,I_3\dots$ be independant indicator variables with a common
parameter $p\in(0,1]$. But here we suppose to have absolutely no
information about the parameter $p$.

Consider an homogeneous Poisson process with rate~1 on $[0,T]$. Let
$N_t$ the number of indicator variables observed up to time $t$. So
$N_t$ counts that number of points in the Poisson process in which an
indicator is observed. We want to stop on the last indicator equal to
1 which arrived in the time interval [0,T].

It is well-known that independent thnning of Poisson process is again
a Poisson process. Hence the arrival process of the 1's is a Poisson
process with unknown rate $p$. Let $\tilde N$ denote the thinned
Poisson process of successes. Here $\tilde N_t$ counts the number of
points in the Poisson process in which an indicator variable of value
1 is observed. We thus want to stop on the last arrival of this
process $\tilde N$ in the interval [0,T]. Here we suppose that only
this process $(\tilde N_t)$ is observable.

We follow the approach of the recent paper of Bruss and Yor~(2012) to
the so-called \emph{last arrival problem} and derive the optimal
strategy.

Noting that $E[\tilde N_t] = pE[N_t] = pt$ we can follow the reasoning
of Bruss and Yor to conclude that, for all $s,t>0$, the process
$(\tilde N_t)$ must satisfy
\begin{equation}
  E[\tilde N_{t+s} - \tilde N_t | \mathcal{F}_t] = sp = \frac{s}{t}pt
  \,, \label{eq:9}
\end{equation}
where $\mathcal{F}_t$ denotes the filtration defined by $\mathcal{F}_t
= \sigma\{\tilde N_t : 0<u\le t\}$. Bruss and Yor~(2012) called such a
process a p.i.-process, that is a process with proportional
increments. Such a process must satisfy
\begin{equation}
  \label{eq:10}
  E(\tilde N_{t+s}-\tilde N_t | \mathcal{F}_t) = \frac{s}{t}\tilde N_t
  \: \text{a.s.}
\end{equation}
According to Theorem~1 from Bruss and Yor~(see page~3243), $(\tilde
N_t/t)$ is a $\mathcal{F}_t$-martingale for $t\ge \tilde T_1$, where
$\tilde T_1$ is the first jump of $\tilde N_t$. (Note that in the
special case where if $\tilde T_1>T$, there is no ``1'' in the
interval [0,T] and we lose by definition.)

Furthermore we note that $(\tilde N_t)$ and $(\tilde N_t/t)$ have exactly
the same jump times. Since $E(\tilde N_t)=pt$ we have $E(\tilde
N_t/t)=p$, and from the martingale property $E(\tilde N_T/T)=E(\tilde
N_t/t)=p$. Exactly as in Bruss and Yor (2012) it is then optimal to
stop at the $k$-th arrival time if and only if
\begin{equation}
  \label{eq:11}
  k \left( \frac{T-\tilde T_k}{\tilde T_k}\right)
  \le 1.
\end{equation}
And if there is no arrival time $\tilde T_k$ in [0,T] for which this
condition is verified such $k$ then we lose. \qed

\medskip\noindent%
\textbf{Remark.} This is, as far as we are aware, the only case when
the unknown odds-problem allow for a solution which is proved to be
optimal. It would be interesting to know whether the
criterion~(\ref{eq:11}) would also be optimal if both $(N_t)$ and
$(\tilde N_t)$ are observable, because in this case the relevant
filtration would be the larger filtration generated by $(N_u)_{u\le
  t}$ and the indicators seen up to time $t$.

\section{Applications}
\label{sec:applications}

We now give other important applications which can be solved by the
odds algorithm or its newer developments.

\subsection{The ballot problem}
\label{sec:ballot-problem}

Tamaki (2001) considers the problem of stopping on the maximum point
of a random trajectory. One of the models he presents is quickly
solved by the odds theorem. The paper studies also several other
problems and contains interesting ideas for future research.

\subsection{Online Calibration in Local Search (Bontempi, 2011)}
\label{sec:optim-calibr-local}

To search the minimum of a real-valued function a computer creates a
grid of points in the domain of the function and evaluates the
function in each of these points. If the function does not fluctuate
too much we will look at the point giving the smallest value
for our function and think that this point should be close to the real
minimum of the function into the considered domain.

When the grid is really tight we evaluate many points and our estimate
for the minimum becomes better. But this is computationnally
unefficient, and one prefers stochastic approaches.

Bontempi suggests to start from an initial best point $x_0$ and try
random points in the neighbourhood of $x_0$. If the function $f$
evaluates smaller in one of these points, say $x_1$, \emph{follow the
  path $x_0\to x_1$} and hope that there is another better point in
the beighbourhood of $x_1$. This is the new current best point (or
solution).

When to stop? It would be best if one could stop searching when being
on the very best solution. If the neighbourhood of some $x_k$ did not
give better points, go back to $x_{k-1}$ and investigate the second
best solution's neighbourhood. This might give a better solution but
also might not. If this \emph{go backward-take next best} procedure
does not give a better solution, then $x_k$ was indeed the best
solution in the neighbourhood of $x_{k-1}$ but also the best of all
$x_1,x_2,\dots,x_k$. So this a good candidate for a minimum of $f$.

\subsection{Automation and Maintenance (Iung, Levrat, Monnin, Thomas
  (2006, 2007, 2008))}
\label{sec:autom-maint}

The research group of Iung, Monnin, Levrat and Thomas of the
\textit{Centre de Recherche en Automatique de Nancy,} CNRS, France,
has applied the odds algorithm to problems of \emph{automation and
  maintenance}. This application is intended to provide a strategy for
a maintenance tool (a robot) to choose which part of a system to
replace if there is more than one failure. The choice would be based
on the life expectation of this piece of the system, the time used to
replace it and the probability of breakdown of it.

Putting the maintenance problem into the \emph{odds framework} is a
most interesting task (see \cite{thomas07:_lalgor_bruss},
\cite{thomas2008maintenance}, \cite{thomas2008overview} and
\cite{thomas06}) with many challenging questions.

The Nancy research group used the odds algorithm to formulate a
strategy to select the priorities of replacements of parts. They were
aware if the fact that the independence condition needed in the
hypotheses is not always satisfied, because failures may be dependent
of each other.

\subsection{Odds and software}
\label{sec:odds-software}

In Skroch and Turowski~(2010), the odds algorithm is used as a
decision tool for optimal selection when a maintenance task must be
performed on particular software systems.

The authors explain that advanced software systems can reconfigure
themselves at run-time by choosing between alternative options for
performing certain functions and that such options can be built into
the systems. However, they point out that these software systems are
also externally available on open and uncontrolled platforms, such as
Web services and mashups on the Internet.
  
The authors show how run-time software self-adaptation with
uncontrolled external options can be optimized by stopping theory,
yielding the best possible lower probability bound for choosing an
optimal option.

\subsection{Investment models (Bruss and Ferguson (2002))}
\label{sec:investment-models}

In a venture capital investment, one wants to invest a certain amount
of money in, for example, a particular domain of technology. This
money often must have been placed before a fixed date, and therefore
it is highly preferable that the best innovation within this period is
the one that was chosen.

It is not too hard to see that again we are here waiting for some
particular event among all observable events, where we try to detect
the last one, when it happens. The event of interest would here be
described as follows: today is an ``opportunity'' if today's
techonology is better than the ones we observe since the beginning of
our observing period. Call this ``opportunity'' a success and a day
without opportunity a ``failure'', or respectively write 1 and 0. This
is almost an ``odds-theorem'' setting.

For other high risk investment models, see also \L ebek and
Szajowski~(2007).

\newpage


\end{document}